\documentclass[letterpaper,final,twoside,reqno]{amsart}
\usepackage{xspace}
\usepackage{amssymb}
\usepackage[all]{xy}

\newcommand{\baseRing}[1]{\ensuremath{\mathbb{#1}}}
\newcommand{\Z}{\baseRing{Z}}
\newcommand{\Q}{\baseRing{Q}}
\newcommand{\N}{\baseRing{N}}
\newcommand{\R}{\baseRing{R}}
\newcommand{\C}{\baseRing{C}}
\newcommand{\CP}{\baseRing{P}}
\newcommand{\jgg}{\ensuremath{\mathfrak{g}}\xspace}

\newcommand{\jdef}[1]{\emph{#1}}
\newcommand{\HO}{\ensuremath{H_{orb}}}

\newcommand{\cupO}{\ensuremath{\smallsmile_{orb}}}
\newcommand{\PR}{\ensuremath{o}}

\newcommand{\Script}[1]{\ensuremath{{\mathcal{#1}}}}
\newcommand{\KK}{\Script{K}}

\newcommand{\ti}[1]{\tilde{#1}}
\newcommand{\conj}{\overline}

\newcommand{\gl}{\ensuremath{\mathfrak{gl}}}

\newcommand{\DD}{\ensuremath{D}\xspace}

\theoremstyle{plain}
\newtheorem{theorem}{Theorem}[section]
\newtheorem{corollary}[theorem]{Corollary}
\newtheorem{proposition}[theorem]{Proposition}

\theoremstyle{definition}
\newtheorem{definition}[theorem]{Definition}
\newtheorem{remark}[theorem]{Remark}
\newtheorem{example}[theorem]{Example}

\numberwithin{equation}{section}

\DeclareMathOperator{\sym}{Sym}
\DeclareMathOperator{\im}{Im}
\DeclareMathOperator{\gr}{Gr}

\CompileMatrices

\begin{document}

\title{Hodge structures for orbifold cohomology}

\author{Javier Fernandez}
\address{Department of Mathematics\\ University of Utah\\ Salt Lake
  City\\ UT 84112--0090\\USA}

\curraddr{Instituto Balseiro\\Univerisdad Nacional de Cuyo --
  C.N.E.A.\\ Bariloche\\ R{\'\i}o Negro\\ R8402AGP \\Rep\'ublica
  Argentina}

\email{jfernand@cab.cnea.gov.ar}



\date{\today}


\bibliographystyle{amsplain}


\begin{abstract}
  We construct a polarized Hodge structure on the primitive part of
  Chen and Ruan's orbifold cohomology $\HO^k(X)$ for projective
  $SL$-orbifolds $X$ satisfying a ``Hard Lefschetz Condition''.
  Furthermore, the total cohomology $\oplus_{p,q}\HO^{p,q}(X)$ forms a
  mixed Hodge structure that is polarized by every element of the
  K\"ahler cone of $X$. Using results of Cattani-Kaplan-Schmid this
  implies the existence of an abstract polarized variation of Hodge
  structure on the complexified K\"ahler cone of $X$.
  
  This construction should be considered as the analogue of the
  abstract polarized variation of Hodge structure that can be attached
  to the singular cohomology of a crepant resolution of $X$, in the
  light of the conjectural correspondence between the (quantum)
  orbifold cohomology and the (quantum) cohomology of a crepant
  resolution.
\end{abstract}

\maketitle



\section{Introduction}
\label{sec:intro}

The cohomology of a projective variety $X$ underlies several very rich
algebraic structures. For example, if $X$ is smooth, the Hodge
decomposition of $H^k(X,\R)$ defines a pure Hodge structure that is,
in fact, polarized when restricted to its primitive part. In general,
a result of P. Deligne constructs a mixed Hodge structure over
$H^k(X,\R)$. Also, the total cohomology of a smooth projective variety
defines a mixed Hodge structure polarized by any K\"ahler class of
$X$. A result of J. Steenbrink and M. Saito proves that this is still
true for projective orbifolds.

Motivated by some ideas originated in physics, W. Chen and Y. Ruan
defined a new cohomology theory for orbifolds: the \jdef{orbifold
  cohomology}. The purpose of this note is to explore under what
conditions the total orbifold cohomology of a projective orbifold has
a Hodge structure similar to that of a smooth projective variety. Our
main result is that if $X$ is a projective $SL$-orbifold satisfying a
certain condition on the ages of its local groups --that we call the
Hard Lefschetz Condition--, then $\HO^*(X,\R)$ underlies a mixed Hodge
structure that is polarized by any K\"ahler class. This result will
follow from the Hard Lefschetz theorem, the Lefschetz decomposition
and the Hodge-Riemann bilinear relations, all of which will be seen to
hold.  As a consequence of this result and the work of E. Cattani, A.
Kaplan and W. Schmid we also conclude that under the previous
conditions there is a natural polarized variation of Hodge structures
defined over the complexified K\"ahler cone of $X$. This variation is
the asymptotic approximation of a conjectural $A$-model variation of
Hodge structure~\cite{bo:CK-mirror,ar:CF-sqm} for Calabi-Yau
orbifolds.

Our results show that there is a close similarity between the Hodge
structures on the orbifold cohomology of a projective orbifold
satisfying the Hard Lefschetz Condition and the standard cohomology of
a smooth projective variety. This similarity can be considered as
indirect evidence for the several conjectures posed by
Ruan~\cite{ar:ruan-cohomology_ring_crepant} regarding the existence of
equivalences between the orbifold cohomology of an orbifold $X$ and
the cohomology of a crepant resolution of $X$.

We will review some notions of Hodge theory and orbifold cohomology in
Sections~\ref{sec:hodge_theory}
and~\ref{sec:orbifold_cohomology}. The Hard Lefschetz Theorem is
discussed in Section~\ref{sec:hard_lefschetz_condition}, while
Section~\ref{sec:hodge_structure} analyzes the polarized Hodge
structure on the orbifold cohomology of appropriate orbifolds.



\section{Hodge theory preliminaries}
\label{sec:hodge_theory}

In this section we recall the basic notions of Hodge theory.
Additionally we state some classical results for the de Rham and
Dolbeault cohomologies of a projective orbifold. We refer
to~\cite{ar:CK-luminy,ar:CKS,bo:griffiths-topics,ar:S-vhs} for more
details and proofs.

Let $V$ be a finite dimensional $\R$-vector space and $k\in\Z$. A
\jdef{(real) Hodge structure of weight $k$} on $V$ is a decomposition
of $V_\C := \C\otimes V$, $V_\C = \oplus_p H^{p,k-p}$ such that
$\conj{H^{p,k-p}} = H^{k-p,p}$ for all $p$. Alternatively, such Hodge
structure is described by a decreasing filtration $F$ of $V_\C$ such
that $V_\C = F^p \oplus \conj{F^{k-p+1}} \text{ for all }p$. The
relation between $F$ and $H$ is that $F^p = \oplus_{a\geq p}
H^{a,k-a}$ while $H^{a,k-a} = F^a \cap \conj{F^{k-a}}$.  The numbers
$h^{p,q} := \dim H^{p,q}$ are called the \jdef{Hodge numbers} of the
structure.

A Hodge structure $(V,H,k)$, is~\jdef{polarized} by the nondegenerate,
bilinear, $(-1)^k$-symmetric form $Q$ on $V$ if
\begin{equation}
  \label{eq:PHS-orthogonality-H}
  Q(H^{p,k-p},H^{q,k-q})\ =\ 0 \quad \text{ unless } p+q=k
\end{equation}
\begin{equation}
  \label{eq:PHS-positivity-H}
  Q(C v, \conj{v}) \ >\ 0 \quad \text{ for all } v\in V_\C-\{0\},
\end{equation}
where $C:V_\C\rightarrow V_\C$ is the \jdef{Weil operator} defined by
$C(v) := i^{p-q} v$ for $v\in H^{p,q}$.

\begin{example}
  For a smooth projective or compact K\"ahler manifold $X$ of
  dimension $n$ with a choice of K\"ahler form $\omega$ and
  $k=0,\ldots, n$, the primitive part of the cohomology $H^k_\PR(X,\C)$
  (see~\eqref{eq:primitive-H}) is a Hodge structure of weight $k$
  polarized by the form
  \begin{equation}
    \label{eq:polarization-Kahler-l}
    Q_{k}(\alpha,\beta)\ :=\ Q(\alpha,\beta \wedge \omega^{n-k}), \quad \text{
      if } \alpha, \beta \in H^k(X,\C)
  \end{equation}
  where 
  \begin{equation}
    \label{eq:polarization_kahler}
    Q(\alpha,\beta)\ :=\ {(-1)}^{k(k-1)/2} \int_X \alpha \wedge \beta,\quad
    \text{ for } \alpha\in H^k(X,\C).
  \end{equation}
\end{example}

The primitive cohomology of a smooth projective variety $X$ is the
typical example of a polarized Hodge structure. The same properties can
be extended to the case where $X$ has some mild singularities. This is
the case when $X$ is a projective orbifold, that is, an orbifold that
can be realized as a projective variety. In what follows we will be
mostly interested in this setup. We refer
to~\cite{ar:satake-orbifold_def,ar:Chen_Ruan-orbifold_cohomology,ar:Chen_Ruan-orbifold_gromov_witten}
for the notion of orbifold.

Let $X$ be a projective orbifold of dimension $n$ with K\"ahler class
$\omega \in H^{1,1}(X,\R) := H^{1,1}(X) \cap H^2(X,\R)$. Define     
\begin{equation*}
  L_\omega: H^*(X,\R) \rightarrow H^*(X,\R), \quad
  L_\omega(\alpha) \ :=\ \omega \wedge \alpha.
\end{equation*}
The \jdef{primitive cohomology} of $X$ is defined as
\begin{equation}
  \label{eq:primitive-H}
  H^p_\PR(X,\C) :=
  \ker(L_\omega^{n-p+1}:H^{p}(X,\C) \rightarrow H^{2n-p+2}(X,\C))
\end{equation}
for $p=0,\ldots, n$ and $\{0\}$ otherwise.

\begin{theorem}\label{th:polarized_hodge_orbifold} 
  Let $X$ be a projective orbifold of dimension $n$ with
  K\"ahler class $\omega\in H^{1,1}(X,\R)$.
  \begin{enumerate}
  \item {\em (Hard Lefschetz)\/} For all $p\in \N$ the map
    $L_\omega^p$ induces an isomorphism between $H^{n-p}(X,\R)$ and
    $H^{n+p}(X,\R)$.
  \item {\em (Lefschetz decomposition)\/} For $k=0,\ldots,n$:
    \begin{equation}
      \label{eq:lefschetz_decomposition}
      H^{k}(X,\C) \ =\ \oplus_{p\geq 0} L_\omega^p H_\PR^{k-2p}(X,\C).
    \end{equation}
  \item For all $k=0,\ldots,2n$ there is a decomposition
    \begin{equation}
      \label{eq:Hodge_decomposition_Kahler}
      H^{k}(X,\C)\ =\ \oplus_{p} H^{p,k-p}(X).
    \end{equation}
    For $k=0,\ldots, n$,~\eqref{eq:Hodge_decomposition_Kahler}
    restricts to the primitive cohomology to give
    \begin{equation}
      \label{eq:primitive_Hodge_decomposition}
      H^{k}_\PR(X,\C)\ =\ \oplus_p H^{p,k-p}_\PR(X),
    \end{equation}
    where $H^{p,k-p}_\PR(X) = H^{p,k-p}(X) \cap H^{k}_\PR(X,\C)$.
    Both decompositions~\eqref{eq:Hodge_decomposition_Kahler}
    and~\eqref{eq:primitive_Hodge_decomposition} define Hodge
    structures of weight $k$ on the underlying $H^{k}(X,\R)$ and
    $H^{k}_\PR(X,\R)$ respectively.
  \item The Hodge structure on the primitive
    cohomology~\eqref{eq:primitive_Hodge_decomposition} is polarized
    by the form~\eqref{eq:polarization-Kahler-l}.
  \end{enumerate}
\end{theorem}

Theorem~\ref{th:polarized_hodge_orbifold} is due to J.
Steenbrink~\cite[\S 1]{ar:steenbrink-vanishing_cohomology}, completed
by M. Saito~\cite{ar:saito-mixed_hodge_modules}.  In
Section~\ref{sec:hodge_structure} we will derive the analogous of
Theorem~\ref{th:polarized_hodge_orbifold} for the orbifold cohomology.

A \jdef{mixed Hodge Structure} (MHS) on $V$ consists of a pair of
filtrations of $V_\C$, $(W, F)$, $W$ defined over $\R$ and increasing,
$F$ decreasing, such that $F$ induces a Hodge structure of weight $p$
on $\gr_p^{W}$ for each $p$.  Equivalently, a MHS on $V$ is a
bigrading $ V_\C = \oplus I^{p,q}$ satisfying $I^{p,q}\equiv
\conj{I^{q,p}} \mod(\oplus_{a<p,b<q} I^{a,b})$
(see~\cite[Theorem~2.13]{ar:CKS}).  Given such a bigrading we define:
$ W_l = \oplus_{p+q \leq l} I^{p,q}$, $F^a = \oplus_{p\geq a}
I^{p,q}$.  A MHS is said to \jdef{split} over $\R$ if $I^{p,q}=
\conj{I^{q,p}}$. A map $T \in \gl(V_\C)$ such that $T(I^{p,q}) \subset
I^{p+a,q+b}$ is called a morphism of bidegree $(a,b)$.

A nilpotent linear transformation $N\in\gl(V)$ defines an increasing
filtration, the \jdef{weight filtration}, $W(N) $ of $V$, uniquely
characterized by requiring that, for all $l$, $N(W_l(N))\subset
W_{l-2}(N)$ and that $N^l:\gr_{l}^{W(N)}\rightarrow \gr_{-l}^{W(N)}$
be an isomorphism.

\begin{definition}
  A \jdef{polarized mixed Hodge structure}
  (PMHS)~\cite[(2.4)]{ar:CK-polarized} of weight $k$ on $V$ consists
  of a MHS $(W,F)$ on $V$, a nondegenerate, bilinear,
  ${(-1)}^k$-symmetric form $Q$, and a $(-1,-1)$-morphism $N\in
  \jgg_\R$, the Lie algebra of $O(V,Q)$, such that
  \begin{enumerate}
  \item $N^{k+1}=0$,
  \item $W = W(N)[-k]$, where ${W[-k]}_j = W_{j-k}$,
  \item $Q(F^a,F^{k-a+1}) = 0$ and,
  \item the Hodge structure of weight $k+l$ induced by $F$ on
    $\ker(N^{l+1}:\gr_{k+l}^{W}\rightarrow \gr_{k-l-2}^{W})$ is
    polarized by $Q(\cdot,N^l \cdot)$.
  \end{enumerate}
\end{definition}

Polarized MHSs arise naturally as limits of polarized variations of
Hodge structure. Conversely, a PMHS generates a nilpotent orbit.
Indeed, there is an equivalence between nilpotent orbits and mixed
Hodge structures polarized by an abelian family of operators as
explained by Theorem (2.3) in~\cite{ar:CK-luminy}.

\begin{example}\label{ex:totalcohomology}
  Let $X$ be an $n$-dimensional, smooth projective variety.  Let $V =
  H^*(X,\R)$.  The bigrading $I^{p,q} := H^{n-q,n-p}(X)$ defines a MHS
  on $V$ which splits over $\R$.  The weight and Hodge filtrations are
  then
  \begin{equation*}
    W_l \ =\ \bigoplus_{d\geq 2n-l} H^d(X,\C),\quad F^p \ =\
    \bigoplus_{r}\bigoplus_{s\leq n-p}H^{r,s}(X).
  \end{equation*}
  Given a K\"ahler class $\omega\in H^{1,1}(X,\R)$, let $L_{\omega}\in
  \gl(V_\R)$ denote multiplication by $\omega$.  Note that
  $L_{\omega}$ is an infinitesimal automorphism of the form $Q$
  defined in~\eqref{eq:polarization_kahler} and is a $(-1,-1)$
  morphism of $(W,F)$.  Moreover, the Hard Lefschetz Theorem and the
  Riemann bilinear relations are equivalent to the assertion that
  $L_\omega$ together with $Q$ polarize $(W,F)$.  Let $\KK \subset
  H^{1,1}(X,\R)$ denote the K\"ahler cone and
  \begin{equation*}
    \KK_\C \ :=\ H^{1,1}(X,\R) \oplus i \KK \subset H^2(X,\C)
  \end{equation*}
  the complexified K\"ahler cone.  It then follows from Theorem~(2.3)
  in~\cite{ar:CK-luminy} that for every $\xi\in \KK_\C$, the
  filtration $\exp(L_\xi)\cdot F$ is a Hodge structure of weight $n$
  on $V$ polarized by $Q$.  The map $\, \xi\in \KK_\C \mapsto
  \exp(L_\xi)\cdot F \,$ is the period map (in fact, the nilpotent
  orbit) of a variation of Hodge structure over $\KK_\C$.
\end{example}

Theorem~\ref{th:PMHS_on_HO} and Corollary~\ref{cor:nilpotent_for_HO}
in Section~\ref{sec:hodge_structure} will prove similar properties for
the orbifold cohomology.



\section{Orbifold Cohomology}
\label{sec:orbifold_cohomology}

In this section we briefly review Chen and Ruan's construction of
orbifold cohomology~\cite{ar:Chen_Ruan-orbifold_cohomology}. 

Recall that to an orbifold $X$ we can associate another orbifold,
known as the \jdef{inertia orbifold}
\begin{equation*}
  \ti{X} \ :=\ \{ (p,(g)) : p \in X, g \in G_p\},
\end{equation*}
where $(g)$ denotes the conjugacy class of $g$ in $G_p$, the local
group of $X$ at $p$. If $\{(V_p,G_p,\pi_p):p\in X\}$ is a uniformizing
system for $X$, a uniformizing system for $\ti{X}$ is $\{(V_p^g,
C(g),\pi_{p,g}): (p,(g)) \in \ti{X} \}$, where $V_p^g$ is the fixed
point set of $g$, and $C(g)$ is the centralizer of $g$ in $G_p$. The
topology on $\ti{X}$ is defined so that the sets
$\pi_{p,g}(V_p^g)\simeq V_p^g/ C(g)$ are open.

Even if $X$ is a connected space, $\ti{X}$ need not be. In general
$\ti{X}$ decomposes in connected components
\begin{equation*}
  \ti{X} \ =\ \bigsqcup_{t \in T} X_{t}
\end{equation*}
where $t$ labels each component and $T$ is the set of all such labels.
The components are orbifolds that are compact (complex) if $X$ is
compact (complex). Also, if $X$ is projective, so are the components
(use~\cite{ar:baily-on_the_imbedding}).

We denote by $T^0\subset T$ the set of connected components that
contain points of the form $(p,(1))$, where $1\in G_p$ is the
identity. Also, $X_1 := \bigsqcup_{t\in T^0} X_t$. $X_1$ is called the
\jdef{non-twisted sector}, while the other connected components of
$\ti{X}$ are the \jdef{twisted sectors}.

Define $\pi: \ti{X}\rightarrow X$ by $\pi(p,(g)):= p$ and let $\pi_t$
denote the restriction of $\pi$ to the sector $X_t$.  Notice that
$\pi_1:X_1\rightarrow X$ is an isomorphism; hence, when $X$ is
connected, so is $X_1$. Also, define the involution $I:\ti{X}
\rightarrow \ti{X}$ by $I(p,(g)) := (p,(g^{-1}))$. We also denote by
$I$ the involution induced on $T$.

For each $g \in G_p$ we consider the action $\rho_p(g)$ induced by $g$
on $T_p X$. The eigenvalues of $\rho_p(g)$ are of the form $\exp(2\pi
i \frac{m_j}{m_g})$, where $m_g$ is the order of $g$ and $m_j\in \Z$,
$0\leq m_j< m_g$. Then, the \jdef{index of $g$} at $p$ (also known as
the \jdef{age} or \jdef{degree shifting number} of $g$) is defined by
\begin{equation*}
  i_{(p,g)} \ :=\ \sum_{j=1}^n \frac{m_j}{m_g} \in\Q \quad \text{ where } \quad 
  n\ :=\ \dim_\C X.
\end{equation*}
This index defines a locally constant function of $p$. If $X_t$ is the
sector containing $(p,(g))$, we define $i_{t} := i_{(p,g)}$.
Then~\cite[Lemma 3.2.1]{ar:Chen_Ruan-orbifold_cohomology} shows that
$i_{(p,g)}$ is integral if and only if $\rho_p(g) \in SL(n,\C)$ and
that the (complex) dimension of $X_{t}$ is
\begin{equation}\label{eq:dims}
  \dim_\C X_{t}\ =\ \dim_\C X - i_{t} - i_{I(t)}.
\end{equation}

An orbifold for which $\rho_p(g) \in SL(n,\C)$ for all $p$ and $g$ is
called an \jdef{$SL$-orbifold}. In algebro-geometric terms, this is to
say that (the geometric space associated to) $X$ is a Gorenstein
variety. In particular, for instance, if $X$ is Calabi-Yau, then $X$
is an $SL$-orbifold.

The \jdef{orbifold cohomology groups} of the orbifold $X$ are defined
by
\begin{equation*}
  \HO^k(X,\R) \ :=\ \oplus_{t\in T} H^{k-2i_{t}}(X_{t},\R)
\end{equation*}
where the cohomology groups on the right are the singular cohomology
groups with real coefficients, which are isomorphic to the
corresponding de Rham cohomology groups~\cite{ar:satake-orbifold_def}.
If $X$ is complex then the \jdef{orbifold Dolbeault cohomology groups}
are defined by
\begin{equation*}
  \HO^{p,q}(X)\ :=\ \oplus_{t\in T} H^{p-i_{t},q-i_{t}}(X_{t}).
\end{equation*}

Chen and Ruan define a product structure on orbifold cohomology.  We
will not go into all the details of this definition for which we
refer, once more, to~\cite{ar:Chen_Ruan-orbifold_cohomology}. The
orbifold product defines a rich structure, as the following result
shows~\cite[Theorems 4.1.5 and
4.1.7]{ar:Chen_Ruan-orbifold_cohomology}.

\begin{theorem}
  Let $X$ be a compact, complex orbifold. Then the orbifold product is
  bigraded
  \begin{equation*}
    \cupO: \HO^{p,q}(X) \times \HO^{p',q'}(X) \rightarrow
    \HO^{p+p',q+q'}(X). 
  \end{equation*}
  The total cohomologies $\HO^*(X)$ and $\HO^{*,*}(X)$ become rings
  under $\cupO$, with unit $e_0^X \in\HO^{0,0}(X)$. When $X$ is an
  $SL$-orbifold both rings are super-commutative. Finally, the
  restriction of $\cupO$ to the non-twisted sector coincides with the
  cup product on $H^*(X_1)$ and $H^{*,*}(X_1)$.
\end{theorem}

The cohomology of the twisted sector $X_t$, $H^*(X_t)$, becomes a
$H^*(X_1)$-module under $\cupO$. The following result relates this
module structure to the standard cup product on $H^*(X_t)$, and will
be useful to study the left (orbifold) multiplication by K\"ahler
classes. 

\begin{proposition}\label{prop:products_are_equal}
  Let $X$ be a compact orbifold, $\alpha\in H^p(X_{1})$ and $\beta\in
  H^q(X_{t})$ for $t\in T$. Then $\alpha \cupO \beta = {(\pi_{1}^{-1}
  \circ \pi_t)}^*(\alpha) \wedge \beta \in H^{p+q}(X_{t})$.
\end{proposition}

\begin{proof}
  In the language of~\cite[\S 4]{ar:Chen_Ruan-orbifold_cohomology},
  the obstruction bundle is trivial for dimensional reasons and the
  proposition follows.
\end{proof}



\section{The Hard Lefschetz condition}
\label{sec:hard_lefschetz_condition}

In this section we want to find a necessary and sufficient condition
for the Hard Lefschetz Theorem to hold on the orbifold cohomology.

As we mentioned in the Introduction, it is generally believed that,
under some yet not well understood conditions, there should be some
kind of equivalence between the orbifold cohomology of an orbifold $X$
and the cohomology of a crepant resolution of $X$. For instance, if
$Y\rightarrow X$ is a hyperk\"ahler resolution, Ruan conjectured such
a
relation~\cite{ar:ruan-stringy_topology_orbifolds,ar:ruan-cohomology_ring_crepant}.
In any case, if such equivalence exists, the algebraic structure of
the cohomology of the resolution should have an analogue in the
orbifold cohomology of $X$. One ``feature'' of the singular cohomology
of the resolution is the Hard Lefschetz Theorem. We will concentrate
here in the analogous result for the orbifold cohomology.

If $X$ is a K\"ahler orbifold, for any K\"ahler class $\omega \in
H^{1,1}(X,\R)$ we define the operator
\begin{equation}
  \label{eq:operator_L}
  L_\omega:\HO^*(X)\rightarrow \HO^*(X), \quad L_\omega(\alpha) \ :=\
  \pi_{1}^* \omega \cupO \alpha.
\end{equation}

As we mentioned in section~\ref{sec:hodge_theory}, the Hard Lefschetz
theorem states that for any K\"ahler class $\omega$ on the projective
orbifold $X$ of complex dimension $n$, for all $p\in \N$ the map
$L_\omega^p$ induces an isomorphism between $\HO^{n-p}(X,\R)$ and
$\HO^{n+p}(X,\R)$. To understand this condition, we remember that the
orbifold product by $\omega^p \in H^{p,p}(X) \simeq H^{p,p}(X_{1})$
preserves the forms on each sector, that is, for all $t\in T$,
$L_\omega^p(H^*(X_{t})) \subset H^*(X_{t})$. Then if $L_\omega^p$ is
an isomorphism, it should induce isomorphisms on each $H^*(X_{t})$.
But care must be taken regarding the degrees: $L_\omega^p$ pairs
$\HO^{n-p}(X,\R)$ and $\HO^{n+p}(X,\R)$. So, for each $t$,
$L_\omega^p$ is pairing $H^{n-p-2i_{t}}(X_{t})$ with
$H^{n+p-2i_{t}}(X_{t})$. It is easy to check using~\eqref{eq:dims}
that this is not possible unless
\begin{equation}
  \label{eq:hard_lefschetz_condition}
  i_{t}\ =\ i_{I(t)} \quad \text{ for all } t\in T.
\end{equation}
Finally, it follows immediately from
Theorem~\ref{th:polarized_hodge_orbifold} that
when~\eqref{eq:hard_lefschetz_condition} holds, $L_\omega^p$
are isomorphisms for all $p$. Thus we have proved

\begin{proposition}\label{prop:hard_lefschetz_orbifold}
  Let $X$ be a projective orbifold of dimension $n$ that satisfies the
  Hard Lefschetz Condition~\eqref{eq:hard_lefschetz_condition}. Then
  for each K\"ahler class $\omega$ on $X$ the operator
  \begin{equation*}
    L_\omega^p: \HO^{n-p}(X,\R)\rightarrow \HO^{n+p}(X,\R) 
  \end{equation*}
  is an isomorphism.
\end{proposition}

\begin{remark}
  Condition~\eqref{eq:hard_lefschetz_condition} had already appeared
  in~\cite[\S 4]{ar:ruan-cohomology_ring_crepant} in connection with
  the possibility of defining a ``hermitian'' inner product to fix a
  signature problem in the Kummer surface example. 
\end{remark}

\begin{remark}
  Condition~\eqref{eq:hard_lefschetz_condition} holds trivially for
  orbifolds of dimension $2$ and for orbifolds whose nontrivial
  local groups are isomorphic to $\Z_2$.
\end{remark}

\begin{example}
  Let $\CP_\Delta$ be the simplicial complete Fano toric variety of
  dimension $n$ associated to the reflexive polytope $\Delta \subset
  \R^n$. If $X$ is a generic anticanonical hypersurface of $P$, then
  $X$ is a reduced Calabi-Yau projective orbifold. M.  Poddar shows
  in~\cite[\S 4.2]{ar:poddar-orbifold_hodge_numbers_CY_hypersurface}
  that the twisted sectors of $X$ are isomorphic to $X\cap
  \overline{O_\tau}$, where $\overline{O_\tau}$ is the closure of the torus
  orbit in $\CP_\Delta$ associated to the cone $\tau$ obtained as the
  cone over a face $F^0$ of $\Delta^0$, the \jdef{polar polytope} of
  $\Delta$, and where $1\leq \dim F^0\leq n-2$ (it is easy to see that
  $\dim \overline{O_\tau} = n - (1+\dim F^0)$). Also, it is shown that
  there is, at least, one twisted sector $X_t$ of the previous form
  with $i_t=1$ for each lattice point in the relative interior of
  $F^0$ (actually the result is more specific, but that is not
  required in our application).

  For example, the weighted projective space $\CP(1,1,2,2,6)$ is
  $\CP_\Delta$ where $\Delta$ is the convex hull of $\{(11,-1,-1,-1),
  (-1,-1,5,-1), (-1,5,-1,-1), (-1,-1,-1,-1),$ $ (-1,-1,-1,1)\}$, a
  reflexive polytope in $\R^4$. The dual polytope $\Delta^0$ is the
  convex hull of $w_0=(-1,-2,-2,-6)$, $w_1=(1,0,0,0)$,
  $w_2=(0,1,0,0)$, $w_3=(0,0,1,0)$ and $w_4=(0,0,0,1)$. Then,
  computations show that the only lattice point in the relative
  interior of a face $F^0$ of $\Delta^0$ with $1\leq \dim F^0\leq 2$
  is $(0,-1,-1,-3)$ and $F^0$ is the convex hull of $w_0$ and $w_1$.
  This shows that the dimension of the only twisted sector is $1$, so
  that~\eqref{eq:dims} implies~\eqref{eq:hard_lefschetz_condition} for
  generic anticanonical hypersurfaces of $\CP(1,1,2,2,6)$.

  Analogously, $\CP(1,1,1,3,3) = \CP_\Delta$ where $\Delta$ is the
  convex hull of $\{(8,-1,-1,-1),$ $(-1,-1,2,-1), (-1,2,-1,-1),
  (-1,-1,-1,-1), (-1,-1,-1,2)\}$, and the dual polytope $\Delta^0$ is
  the convex hull of $w_0=(-1,-1,-3,-3)$, $w_1=(1,0,0,0)$,
  $w_2=(0,1,0,0)$, $w_3=(0,0,1,0)$ and $w_4=(0,0,0,1)$. In this case,
  the only lattice point in the relative interior of a face $F^0$ of
  $\Delta^0$ with $1\leq \dim F^0\leq 2$ is $(0,0,-1,-1)$ and $F^0$ is
  the convex hull of $w_0$, $w_1$ and $w_2$. Hence, if $X$ is a
  generic anticanonical hypersurface of $\CP(1,1,1,3,3)$, it has a
  twisted sector $X_{t}$ isomorphic to $X \cap \overline{O_\tau}$,
  with $\dim \overline{O_\tau} = 4 - (1 + \dim F^0) = 1$, so $\dim
  X_{t} = 0$, and, by~\eqref{eq:dims}, $i_t\neq i_{I(t)}$, so
  that~\eqref{eq:hard_lefschetz_condition} fails.

  The same techniques can be applied further to
  test~\eqref{eq:hard_lefschetz_condition} in this toric setting. It
  would be very interesting to understand how common
  condition~\eqref{eq:hard_lefschetz_condition} is.
\end{example}



\section{Hodge structure}
\label{sec:hodge_structure}

In this section we will, first, construct polarized Hodge structures on
the primitive part of $\HO^k(X)$ for appropriate orbifolds. Then we
will see that for these orbifolds the full cohomology carries the
structure described in Example~\ref{ex:totalcohomology}. Finally, we
use this last property to construct a polarized variation of Hodge
structure associated to $X$.

In order to define a Hodge structure on $\HO^k(X)$ a first requirement
is that the cohomology be integrally graded, which immediately leads
to $X$ being an $SL$-orbifold. Also, we need $X$ to be projective
satisfying condition~\eqref{eq:hard_lefschetz_condition}, all
conditions that we will assume in this section.

Let $X$ be a projective orbifold of dimension $n$ and $\omega\in
H^{1,1}(X,\R)$ a K\"ahler class. The \jdef{primitive orbifold
  cohomology} is defined by
\begin{equation*}
  (\HO^p)_\PR(X,\C)\ :=\ \ker(L_\omega^{n-p+1}:\HO^p(X,\C)\rightarrow
  \HO^{2n-p+2}(X,\C) )
\end{equation*}
where $L_\omega$ is defined by~\eqref{eq:operator_L}.  For the
Dolbeault groups we define
\begin{equation*}
  (\HO^{p,q})_\PR(X) \ :=\ \HO^{p,q}(X)\cap (\HO^{p+q})_\PR(X,\C).
\end{equation*}

\begin{theorem}
  \label{prop:PHS_orbifold}
  Let $X$ be a projective $SL$-orbifold of dimension $n$ satisfying
  condition~\eqref{eq:hard_lefschetz_condition}. Then, for each $k$
  the decomposition
  \begin{equation}
    \label{eq:HO_hodge}
    \HO^k(X,\C)\ =\ \oplus_p \HO^{p,k-p}(X)
  \end{equation}
  is a Hodge structure of weight $k$ on $\HO^k(X,\R)$. Furthermore,
  for any K\"ahler class $\omega$, the decomposition
  \begin{equation}
    \label{eq:HO_hodge_primitive}
    (\HO^k)_\PR(X,\C)\ =\ \oplus_p (\HO^{p,k-p})_\PR(X)
  \end{equation}
  defines a Hodge structure of weight $k$ that is polarized by the
  form $Q_k$, defined as the direct sum over all $t\in T$ of
  \begin{equation}
    \label{eq:polarization_pieces}
    Q^{t}_{k-2i_{t}}(\alpha,\beta) \ :=\ 
    \begin{cases}
      Q^{t}(\alpha,\beta\wedge \pi_{t}^*\omega^{n-k}) \text{ if }
      \begin{split}
        & \alpha\in H^{p,k-2i_{t}-p}(X_{t}) \\
        & \beta \in H^{k-2i_{t}-p,p}(X_{t})
      \end{split}\\
      0 \text{ otherwise},
    \end{cases}
  \end{equation}
  where 
  \begin{equation}
    \label{eq:polarization_sectors}
    Q^{t}(\alpha,\beta)\ :=\ (-1)^{k(k-1)/2+i_{t}} \int_{X_{t}} \alpha
    \wedge \beta \quad \text{ for } \alpha \in H^{k-2i_{t}}(X_{t}).
  \end{equation}
\end{theorem}

\begin{proof}
  Since $\omega$ is a K\"ahler class on $X$, $\pi_{t}^*\omega$ is a
  K\"ahler class on the sector $X_{t}$.  By
  Theorem~\ref{th:polarized_hodge_orbifold}, since each $X_t$ is a
  projective orbifold, each $H^{k-2i_{t}}_\PR(X_{t})$ has a Hodge
  structure of weight $k-2i_{t}$ given by
  \begin{equation*}
    H^{k-2i_{t}}_\PR(X_{t},\C)\ =\ \oplus_p
    (H^{p,k-2i_{t}-p}(X_{t}) \cap H^{k-2i_{t}}_\PR(X_{t})),
  \end{equation*}
  and this structure is polarized by the
  form~\eqref{eq:polarization_pieces}.
  
  Tensoring the polarized Hodge structure defined above with
  $\R(i_{t})$ (where $\R(k)$ is $\R$ viewed as a polarized Hodge
    structure of weight $2k$) we obtain a new polarized Hodge
  structure of weight $k$ on $H^{k-2i_{t}}_\PR(X_{t})$:
  \begin{equation*}
    H^{k-2i_{t}}_\PR(X_{t},\C)\ =\ \oplus_p
    (H^{p-i_{t},k-p-i_{t}}(X_{t}) \cap
    H^{k-2i_{t}}_\PR(X_{t},\C)).  
  \end{equation*}
  This structure is still polarized by the same form
  $Q^{t}_{k-2i_{t}}$. 
  
  Taking the direct sum over all sectors $t\in T$ the result follows.
\end{proof}

\begin{remark}
  The Hodge structure on $\HO^k(X)$ is already implicit
  in~\cite[Proposition 3.3.2]{ar:Chen_Ruan-orbifold_cohomology}. An
  important point is that in order to define the polarization we are
  not using the Kodaira-Serre duality used in that paper. Instead, we
  are using integration on $X_{t}$ which is only possible under
  condition~\eqref{eq:hard_lefschetz_condition}. This form, that we
  found naturally in the Hodge theoretic setting, was proposed
  in~\cite{ar:ruan-cohomology_ring_crepant} in connection with the
  signature of the Poincar\'e form for the Kummer surface.
\end{remark}

In Theorem~\ref{prop:PHS_orbifold} we considered the forms $Q^{t}$
defined by~\eqref{eq:polarization_sectors}. We now define $Q$ to be
the direct sum of the forms $Q^{t}$ over all $t\in T$.  A routine
check shows that this form is bilinear, nondegenerate and
$(-1)^n$-symmetric.

\begin{theorem}
  \label{th:PMHS_on_HO}
  Let $X$ be a projective $SL$-orbifold of dimension $n$ satisfying
  the Hard Lefschetz Condition~\eqref{eq:hard_lefschetz_condition} and
  $\omega$ a K\"ahler class. Then the mixed Hodge structure $(W,F)$
  defined by the bigrading $I^{p,q} := \HO^{n-q,n-p}(X)$:
  \begin{equation}
    \label{eq:filtrations_from_bigrading}
    W_l\ :=\ \oplus_{k\geq2n-l} \HO^k(X,\C), \quad F^p\ :=\ \oplus_a
    \oplus_{b\leq n-p} \HO^{a,b}(X)
  \end{equation}
  is a weight $n$ MHS polarized by the operator $L_\omega$ defined
  in~\eqref{eq:operator_L} with the bilinear form $Q$ defined above.
\end{theorem}

\begin{proof}
  Since $\cupO$ preserves the bigrading of the Dolbeault groups,
  $L_\omega$ is a $(-1,-1)$ morphism of $(W,F)$ and, in particular,
  $L_\omega^{n+1}=0$. To check that $L_\omega \in \jgg_\R$ it suffices
  to check that
  \begin{equation*}
    Q^{t}(L_\omega\alpha,\beta) + Q^{t}( \alpha,
    L_\omega\beta)\ =\ 0 \text{ for } \alpha,\beta \in
    H^*(X_{t}) 
  \end{equation*}
  since the decomposition in sectors is $Q$-orthogonal, and this is a
  straightforward computation. The reality of $L_\omega$ follows from
  $\omega$ being real and the reality of $\cupO$.

  To show that $W$ defined by~\eqref{eq:filtrations_from_bigrading} is
  $W(L_\omega)[-n]$ it suffices to check that $L_\omega(W_l)\subset
  W_{l-2}$ for all $l$ and that $L_\omega^j: W_{j+n}/W_{j+n-1}
  \rightarrow W_{-j+n}/W_{-j+n-1}$ are isomorphisms for all $j$. The
  first condition follows immediately from the definitions. The second
  condition is equivalent to requiring that $L_\omega^j :
  \HO^{n-j}(X)\rightarrow \HO^{n+j}(X)$ be isomorphisms, which is true
  by Proposition~\ref{prop:hard_lefschetz_orbifold}.
  
  Notice that, for each $t\in T$, $Q^{t}(H^{a-i_{t},b-i_{t}}(X_{t}) ,
  H^{c-i_{t},d-i_{t}}(X_{t})) = 0$ unless $a+c=n$ which, by the $Q$
  orthogonality of the sector decomposition, implies that
  $Q(\HO^{a,b}(X), \HO^{c,d}(X)) = 0$ if $a+c\leq n-1$. In turn, this
  last condition implies that $Q(F^{p},F^{n-p+1})=0$ for all $p$.
  
  Finally we have to prove that the Hodge structure of weight $n+l$
  induced by $F$ on $\ker(L_\omega^{l+1}:\gr_{n+l}^{W}\rightarrow
  \gr_{n-l-2}^{W})$ is polarized by $Q(\cdot,L_\omega^l
  \cdot)$. Using~\eqref{eq:filtrations_from_bigrading} we are studying
  the Hodge structure on $(\HO^{n-l})_\PR(X)$ given
  by~\eqref{eq:HO_hodge_primitive} (with $k=n-l$) and the
  polarization is precisely the one that, because of
  Theorem~\ref{prop:PHS_orbifold}, polarizes this structure.
\end{proof}

For $V=\oplus_{p=0}^n V_{2p}$ a finite dimensional, graded, $\C$-vector
space with a $\sym V_2$-module structure, a notion of Frobenius module
is introduced in~\cite{ar:CF-sqm}. A natural generalization exists if
$V$ is a $\sym W$-module for a subspace $W\subset V_2$. In this
context, Theorem~\ref{th:PMHS_on_HO} implies that $V:=\oplus_{p=0}^n
\HO^{p,p}(X)$ is a polarizable $\sym H^{1,1}(X)$-Frobenius module of
weight $n:=\dim_\C X$.

An immediate consequence of Theorems~\ref{th:PMHS_on_HO}
and~(2.3) in~\cite{ar:CK-luminy} is the following

\begin{corollary}
  \label{cor:nilpotent_for_HO}
  Let $X$ be a projective $SL$-orbifold of dimension $n$ satisfying
  condition~\eqref{eq:hard_lefschetz_condition}. Given
  $\omega_1,\ldots, \omega_r$ in the K\"ahler cone of $X$ such that
  they form a basis of $H^{1,1}(X)$ and $F$, $Q$ as in
  Theorem~\ref{th:PMHS_on_HO}, then
  \begin{equation*}
    \theta: U^r \rightarrow \DD, \quad \theta(z_1,\ldots,z_r)\ :=\ 
    \exp(\sum_{j=1}^r z_j L_{\omega_j}) \cdot F
  \end{equation*}
  is a nilpotent orbit, where $U$ is the upper half plane $\im z>0$.
  Moreover, if ${\mathcal K}_\C$ is the complexified K\"ahler cone of
  $X$, $\xi \mapsto \exp(L_\xi)\cdot F$ is a polarized variation of
  Hodge structure of weight $n$ defined over ${\mathcal K}_\C$.
\end{corollary}

\begin{proof}
  Since $\cupO$ is commutative for multiplication by $2$-classes,
  $\{L_{\omega_1}, \ldots, L_{\omega_r}\}$ is a commuting subset of
  $\jgg_\R$. By Theorem~\ref{th:PMHS_on_HO} the weight filtration
  $W(\sum_{j=1}^r \lambda_j L_{\omega_j})$ is the constant $W$ defined
  by~\eqref{eq:filtrations_from_bigrading}, hence independent of
  $\lambda_j>0$. Again by the same result, $(W,F)$ is polarized by
  every $\sum_{j=1}^r \lambda_j L_{\omega_j}$ with $\lambda_j>0$.
  Noticing that the mixed Hodge structure $(W,F)$ is split over $\R$,
  Theorem~(2.3) in~\cite{ar:CK-luminy} completes the proof.
\end{proof}





\def\cprime{$'$}
\providecommand{\bysame}{\leavevmode\hbox to3em{\hrulefill}\thinspace}
\providecommand{\MR}{\relax\ifhmode\unskip\space\fi MR }
\providecommand{\MRhref}[2]{%
  \href{http://www.ams.org/mathscinet-getitem?mr=#1}{#2}
}
\providecommand{\href}[2]{#2}


\end{document}